\documentclass[11pt]{article}
\usepackage{url}
\usepackage{graphics}
\usepackage{amsmath}
\usepackage{amsfonts}
\usepackage{amsthm}

\newcommand{\Cr}{{\hbox{cr}}}
\newcommand{\rcr}{\overline{\hbox{cr}}}
\newcommand{\dd}{{\cal D}}
\newcommand{\pid}{{\pi_\dd}}
\newcommand{\thmatrix}{{Q}}

\newcommand\floor[1]{{\lfloor{#1}\rfloor}}
\newcommand{\bigfloor}[1]{{\biggl\lfloor{#1}\biggr\rfloor}}
\newcommand{\smallfloor}[1]{{\bigl\lfloor{#1}\bigr\rfloor}}
\newcommand{\ksn}{{K_{7,n}}}

\newcommand{\LR}{\mathbb{R}}

\newcommand{\cA}{{\cal A}}
\newcommand{\cBB}{{\cal B}}
\newcommand{\cS}{{\cal S}}
\newcommand{\cC}{{\cal C}}

\newcommand{\cN}{{\cal N}}
\newcommand{\cK}{{\cal K}}
\newcommand{\Tr}{\mathrm{Tr}}
\newcommand{\e}{\epsilon}
\newcommand{\hlf}{\frac 12}
\newcommand{\beq}{\begin{equation}}
\newcommand{\eeq}{\end{equation}}
\newcommand{\beann}{\begin{eqnarray*}}
\newcommand{\eeann}{\end{eqnarray*}}
\newcommand{\bc}{\begin{center}}
\newcommand{\ec}{\end{center}}
\newcommand{\reff}[1]{(\ref{#1})}

\newtheorem{theorem}{Theorem}
\newtheorem{lemma}[theorem]{Lemma}

\newtheorem{proposition}[theorem]{Proposition}

\title{Improved bounds for the crossing numbers  of $K_{m,n}$ and $K_n$}
\author{E.~de Klerk\thanks{Department of Combinatorics and Optimization,
Faculty of Mathematics, University of Waterloo,
Waterloo, Canada, N2L 3G1}\and
J.~Maharry\thanks{Department of Mathematics, The Ohio
State University, Columbus OH 43220, USA}\and
D.V.~Pasechnik\thanks{Theoretische Informatik,
FB20 Informatik, J.W.~Goethe-Universit\"at, Robert-Mayer Str.11-15,
Postfach 11 19 32, 60054 Frankfurt(Main),
Germany. Partially supported by the DFG Grant
SCHN-503/2-1. A part of the research completed while
supported by the Mathematical Sciences Research Institute
(MSRI) at Berkeley CA, USA}\and
R.B.~Richter\footnotemark[1]
\and
G.~Salazar\thanks{Instituto de Investigacion en Comunicacion Optica,
Universidad Autonoma de San Luis Potosi,
Av. Karakorum 1470, Lomas 4ta Seccion,
San Luis Potosi, SLP Mexico 78210.
Supported by grants CONACYT J32168 and FAI-UASLP.
A part of the research completed while
on sabbatical leave at The Ohio
State University, Columbus OH, USA}}
\begin{document}
\maketitle
\vspace{-1cm}
\begin{abstract}
It has been long--conjectured that the crossing number $\Cr(K_{m,n})$ of the
complete bipartite graph $K_{m,n}$ equals the Zarankiewicz Number $Z(m,n):=
\floor{\frac{m-1}{2}} \floor{\frac{m}{2}} \floor{\frac{n-1}{2}}
\floor{\frac{n}{2}}$.  Another long--standing conjecture states that the
crossing number $\Cr(K_n)$ of the complete graph $K_n$ equals
$Z(n):=\frac{1}{4}\smallfloor{\frac{n}{2}}
\smallfloor{\frac{n-1}{2}}
\smallfloor{\frac{n-2}{2}}\smallfloor{\frac{n-3}{2}}$.
In this paper we show the following improved bounds on the asymptotic ratios of
these crossing numbers and their conjectured values:
\begin{itemize}
\item[(i)] for each fixed $m\ge
9$, $\lim_{n\to\infty} \Cr(K_{m,n})/Z(m,n) \ge 0.83m/(m-1)$;
\item[(ii)]
$\lim_{n\to\infty} \Cr(K_{n,n})/Z(n,n) \ge 0.83$; and
\item[(iii)] $\lim_{n\to\infty}
\Cr(K_{n})/Z(n) \ge 0.83$.
\end{itemize}
The previous best known lower bounds were
$0.8m/(m-1), 0.8$, and $0.8$, respectively.  These improved bounds are obtained
as a consequence of the new bound $\Cr(\ksn) \ge 2.1796n^2 - 4.5n$.  To obtain
this improved lower bound for $\Cr(\ksn)$, we use some elementary topological
facts on drawings of $K_{2,7}$ to set up a quadratic program on $6!$ variables
whose minimum $p$ satisfies $\Cr(\ksn) \ge (p/2)n^2 - 4.5n$, and then use
state--of--the--art quadratic optimization techniques
combined with a bit of invariant theory of permutation groups
to show that $p \ge 4.3593$.
\end{abstract}


\begin{section}{Introduction}
In the earliest known instance of a crossing number question, Paul Tur\'an
raised the problem of calculating the crossing number  of the complete
bipartite graphs $K_{m,n}$.  
Tur\'an's interesting account of the origin of this problem  
can be found in~\cite{turanwelcome}.

We recall that in a {\em drawing} of a graph in the plane, different vertices
are drawn as different points, and each edge is drawn as a simple arc whose
endpoints coincide with the drawings of the endvertices of the edge.  
Furthermore, the interior of the arc for an edge is disjoint from all the
 vertex points.
We often make no distinction between a graph object, such as a
vertex, edge, or cycle, and the subset of the plane  that represents it in a
drawing of the graph.

The {\em crossing number} $\Cr(G)$ of a graph $G$ is the minimum number of
pairwise intersections of edges (at a point other than a vertex)  in a drawing
of $G$ in the plane.

Exact crossing numbers of graphs are in general very difficult to compute.
Long--standing conjectures involve the crossing numbers of interesting families
of graphs, such as  $K_{m,n}$ and $K_n$.
On a positive note, it was recently proved by 
Glebskii and Salazar~\cite{glebskysalazar} that the
crossing number of the Cartesian product $C_m\times C_n$  of the cycles of
sizes $m$ and $n$ equals its long--conjectured value, namely $(m-2)n$, at least
for $n \ge m(m+1)$.   For a recent survey of crossing number results, see for
instance Shahrokhi {\em et al.}~\cite{sssv}.

Zarankiewicz published a paper \cite{zarankiewicz}, 
in which he claimed that $\Cr(K_{m,n})=Z(m,n)$
for all positive integers $m$, $n$, where
\begin{equation}\label{eq:zmn}
Z(m,n)=
\bigfloor{\frac{m-1}{2}} \bigfloor{\frac{m}{2}}
 \bigfloor{\frac{n-1}{2}} \bigfloor{\frac{n}{2}}.
\end{equation}
However, several years later Ringel and Kainen independently found a hiatus in
Zaran\-kiewicz's argument.  A comprehensive account of the history of the
problem, including a discussion of the gap in Zarankiewicz's argument, is given
by Guy~\cite{declineandfall}.

Figure~\ref{fig:fig1}  shows a drawing of $K_{4,5}$ with $8$ crossings.
As Zarankiewicz observed, such a drawing strategy can be naturally generalized
to construct, for any positive integers $m,n$, drawings of $K_{m,n}$ with
exactly $Z(m,n)$ crossings.
This observation implies the following
well--known upper bound for $\Cr(K_{m,n})$:
\begin{equation*}\label{zarineq}
\Cr(K_{m,n}) \le Z(m,n).
\end{equation*}
\begin{figure}\label{fig:fig1}
\begin{center}
  \resizebox{60mm}{!}{\includegraphics{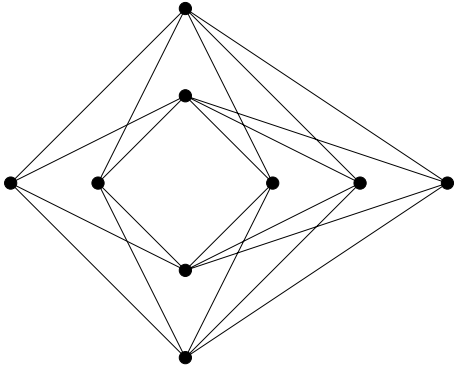}}
  \caption{A drawing of $K_{4,5}$ with $8$ crossings.  A similar strategy can
  be used to construct drawings of $K_{m,n}$ with exactly $Z(m,n)$ crossings.}
\end{center}
\end{figure}
No one has yet exhibited a drawing of any $K_{m,n}$ with fewer than $Z(m,n)$
crossings.  In
allusion to Zarankiewicz's failed attempt to prove that this is the crossing
number of $K_{m,n}$, the following is commonly known as  {\em Zarankiewicz's
Crossing--Number Conjecture}:
\begin{equation*}\label{zarconj}
\Cr(K_{m,n}) \stackrel{?}{=}
Z(m,n),\qquad
%\bigfloor{\frac{m-1}{2}} \bigfloor{\frac{m}{2}}
% \bigfloor{\frac{n-1}{2}} \bigfloor{\frac{n}{2}},
 \text{for all positive integers } m,n.
\end{equation*}

In 1973, Guy and Erd{\H o}s wrote ``[Al]most all questions
that one can ask about crossing numbers remain unsolved'' \cite{erdosguy}.
More than three decades later, despite some definite progress in our
understanding of this elusive parameter, most of the fundamental and more
important questions about crossing numbers remain open.  Zarankiewicz's
Conjecture has been verified by Kleitman~\cite{kleitman}
for $\min\{m,n\} \le 6$ and by Woodall~\cite{woodall} 
for the special cases $7\leq m\leq 8$, $7 \leq n\leq 10$.

Since the crossing number of $K_{m,n}$ is unknown for all other values of $m$
and $n$, it is natural to ask what are the best general lower bounds known for
$\Cr(K_{m,n})$.  A standard counting argument, together with the fact that
$\Cr(K_{5,n})$  is as conjectured, yields the best general lower bound
(\ref{eq:cr-kmn}) known
for $\Cr(K_{m,n})$. It goes as follows: suppose we know
a lower bound $c_r$ on $\Cr(K_{r,n})$, for $2<r<m\leq n$.
Each crossing in the embedding of $K_{m,n}$ lies in $\binom{m-2}{r-2}$
distinct $K_{r,n}\subset K_{m,n}$.
As there are in total $\binom{m}{r}$ distinct $K_{r,n}$'s, one obtains
\begin{equation}\label{eq:cr-kmn}
\Cr(K_{m,n})\geq \frac{c_r\binom{m}{r}}{\binom{m-2}{r-2}},
\quad\text{for $r=5$ one derives}\quad \Cr(K_{m,n}) \ge 0.8\, Z(m,n).
\end{equation}
A small improvement on the $0.8$ factor (roughly to something around
$0.8001$) was recently reported by Nahas~\cite{nahas}.

Zarankiewicz's Conjecture for $K_{7,n}$ states that
\begin{equation}
\Cr(K_{7,n}) \stackrel{?}{=} 9\bigfloor{\frac{n-1}{2}}
\bigfloor{\frac{n}{2}} =
 \left\{ \begin{array}{lrr}
 2.25n^2 -4.5n + 2.25, & \text{$n$ odd,} & n \ge 7 \\
2.25n^2 -4.5n, &  \text{$n$ even,} &  n \ge 8
        \end{array} \right.
\nonumber
\end{equation}

As we observed above, this has been verified only for $n = 7, 8, 9$, and $10$.
Using $\Cr(K_{7,10})= 180$, a standard counting argument gives the best known
lower bounds for $\Cr(K_{7,n})$ for $11 \le n \le 22$.  However, for $n \ge
23$, the best known lower bounds for $\Cr(K_{7,n})$ are obtained by  the same
counting argument, but using the known value of $\Cr(K_{5,n})$ instead of
$\Cr(K_{7,10})$.  Summarizing, previous to this paper, the best known lower
bounds for $\Cr(7,n)$ were:

\begin{equation}\label{previousbound}
\Cr(K_{7,n}) \ge  \left\{ \begin{array}{lr}
          2n(n-1), &  \text{$11\le n \le 22$}, \\
          2.1n^2 -4.2n + 2.1,  &\text{odd }   n \ge 23,\\
          2.1n^2 -4.2n, &  \text{even } n \ge 24,
        \end{array} \right.
\end{equation}

In this paper we prove the following.
\begin{theorem}\label{maintheorem}
For all integers $n$,
\begin{equation}
\Cr(K_{7,n}) > 2.1796 n^2 - 4.5n. \nonumber
\end{equation}
\end{theorem}

An elementary calculation shows that this is an improvement, for all $n \ge
23$, on the bounds for $\Cr(\ksn)$ given in (\ref{previousbound}).

The strategy of the proof can be briefly outlined as follows.  Let $(A,B)$ be
the
bipartition of the vertex set of $\ksn$, where $|A|  = 7$ and $|B| = n \ge 2$.
Let $b, b'$ be vertices in $B$.  In any drawing $\dd$ of $\ksn$, the number of
crossings that involve an edge incident with  $b$ and an edge incident with
$b'$ is bounded from below by a function of the cyclic rotation schemes of $b$
and $b'$.  This elementary topological observation on drawings of $K_{2,7}$
naturally yields a standard quadratic (minimization) program whose minimum $p$
satisfies $\Cr(\ksn) \ge (p/2)n^2 - 4.5n $ (see Lemma~\ref{7ninequality}).
We then use state--of--the--art quadratic programming techniques to show that
$p \ge 4.3593$ (see Lemma~\ref{thebound}), thus implying
Theorem~\ref{maintheorem}.

The rest of this paper is organized as follows.  In Section~\ref{setup}, we
review some elementary topological observations about drawings of $K_{2,n}$,
and use these facts to set up the quadratic program mentioned in
the previous paragraph.   The bound for $\Cr(K_{7,n})$ in terms of the minimum
of this quadratic program is the content of Lemma~\ref{7ninequality}. In
Section~\ref{quadraticprogramming} we prove Proposition~\ref{thebound},
which gives a
lower bound for the quadratic program.   As we observe at the end of
Section~\ref{thebound}, Theorem~\ref{maintheorem} is an obvious consequence of
Lemmas~\ref{7ninequality} and~\ref{thebound}.  In Section~\ref{consequences} we
discuss consequences of Theorem~\ref{maintheorem}: 
the improved bound for $\Cr(\ksn)$ implies improved
asymptotic bounds for the crossing numbers of $\Cr(K_{m,n})$ and $\Cr(K_n)$.
\end{section}

\begin{section}{Quadratic optimization problem yielding a lower
bound for $\Cr(K_{m,n})$}\label{setup}

Our goal in this section is to establish Lemma~\ref{7ninequality}, a statement
that gives a lower bound for $\Cr(K_{m,n})$, for $m\leq n$,
(and thus for $\Cr(\ksn)$) in terms of the solution of a
quadratic minimization problem on $(m-1)!$ variables.

Let $n \ge m$ be fixed.  Let $V$ denote the vertex set of $K_{m,n}$, and let
$(A,B)$ denote the bipartition of $V$ such that each vertex of $A=\{ a_0, a_1,
\ldots, a_{m-1}\}$ is adjacent to each vertex of $B = \{b_0, b_1, \ldots,
b_{n-1}\}$.

Consider a fixed drawing $\dd$ of $K_{m,n}$.
To each vertex $b_i$ we associate a
cyclic ordering $\pid(b_i)$ of the elements in $A$, defined by the (clockwise)
cyclic order  in which the edges incident with $b_i$ leave $b_i$ towards the
vertices in $A$ (see Figure~\ref{fig:fig2}).   Let $\Pi$ denote the set of
all cyclic orderings of $\{a_0, a_1,\ldots, a_{m-1}\}$.  Note that
$|\Pi| = m!/m = (m-1)!$.

Following Kleitman~\cite{kleitman}, let $\Cr_\dd(b_i, b_j)$
denote the number of crossings in $\dd$ that involve an edge incident with
$b_i$ and an edge incident with $b_j$. Further, 
let $\rho_1,\rho_2\in\Pi$ and  
$\thmatrix{(\rho_1,\rho_2)}$ be the minimum number of
interchanges of adjacent elements of $\rho_1$ required to produce
$\rho_2^{-1}$.  Then, for all $b_i,b_j$ with $b_i\ne b_j$,
\begin{equation}\label{bound1} \Cr_\dd(b_i, b_j) \ge
\thmatrix({\pid(b_i),\pid(b_j)}). %, \quad \text{for all } b_i, b_j \text{
%such that } b_i\neq b_j.
\end{equation}
This inequality is stated in \cite{kleitman} and proved in \cite{woodall}.
\begin{figure}
\begin{center}\label{fig:fig2}
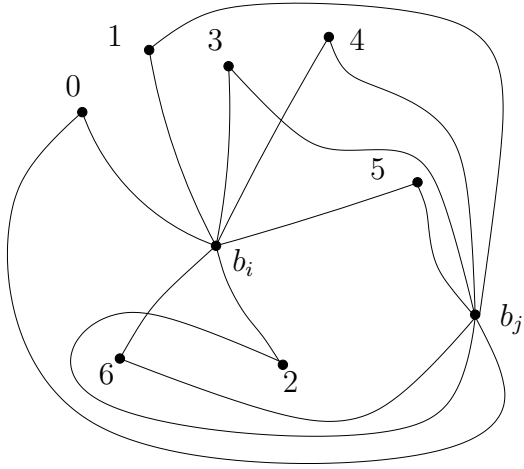
  \caption{Here $m=7$.
Vertices $b_i$ and $b_j$ have cyclic orderings
$(0\, 1\, 3\, 4\, 5\, 2\, 6 )$ and $(0\, 2\, 6\, 5\, 3\, 4\, 1)$,
respectively (we write $i$ for $a_i$ for the sake of brevity).
It is easy to check that the minimum number of interchanges among
adjacent elements in $(0\, 1\, 3\, 4\, 5\, 2\, 6 )$ required to
obtain  $(0\, 2\, 6\, 5\, 3\, 4\, 1)^{-1}$
(namely $(0\, 1\, 4\, 3\, 5\, 6\, 2)$) is $2$.  Thus,
$Q((0\, 1\, 3\, 4\, 5\, 2\, 6 ),(0\, 1\, 3\, 4\, 5\, 2\, 6 )) = 2$.
Therefore, there must be at least two crossings
(as it is indeed the case in the drawing above)
that involve edges incident with $b_i$ and $b_j$.}
\end{center}
\end{figure}
This observation alone yields a lower bound for $\Cr(K_{m,n})$, as follows.
Fix any drawing $\dd$ of $K_{m,n}$.  For each $\rho\in\Pi$, let
\begin{equation*}
x_\rho := \frac{1}{n}{|\{b_i \in B \mid \pi_\dd(b_i) = \rho\}|}
\end{equation*}
The matrix $\thmatrix$ can be viewed
as the matrix of quadratic form $\thmatrix(\cdot,\cdot)$
on the space $\LR^{|\Pi|}$.
It follows from (\ref{bound1}) that
\begin{align*}
{\Cr(\dd)}{}  & \ge     \sum_{\substack{\rho,\rho' \in \Pi \\ \rho\neq \rho'}}
 \thmatrix(\rho, \rho') (x_{\rho}n)(x_{\rho'}n) +  \sum_{\rho\in\Pi}
 \thmatrix(\rho, \rho) \binom{x_{\rho} n}{2}=  \\
%& =  {\frac{n^2}{2}} \sum_{\rho,\rho' \in \Pi}
% \thmatrix(\rho, \rho') x_{\rho} x_{\rho'}   -
% {\frac{n}{2}} \sum_{\rho\in\Pi} \thmatrix(\rho, \rho) {x_{\rho}} =
% \\
 & =  \frac{n}{2}\left( n\sum_{\rho,\rho' \in \Pi}
 \thmatrix(\rho, \rho') x_{\rho} x_{\rho'}  -
 \bigfloor{\frac{m}{2}}\bigfloor{\frac{m-1}{2}}\right),
 \end{align*}
where we have used the (easily verifiable, see e.g.
\cite{woodall}) fact that
$\thmatrix(\rho,\rho) = \floor{m/2} \floor{(m-1)/2}$ for every $\rho \in \Pi$.

Since the drawing $\dd$ was arbitrary, we have proved the following.
\begin{lemma}\label{7ninequality}
Let $Q$ be the $(m-1)!\times (m-1)!$ matrix of the form
$Q(\cdot,\cdot)$,  and let $e$ denote the all ones vector. Then, for every
integer $n \ge m\ge 2$,
\begin{align*}
\Cr(K_{m,n}) & \ge \frac{n}{2} \left(n\min
\{x^T Q x \ \mid \ x \in \mathbb{R}_+^{(m-1)!} ,
e^T x = 1\} -
\bigfloor{\frac{m}{2}}\bigfloor{\frac{m-1}{2}}\right),
\\
\Cr(K_{7,n}) & \ge \frac{n}{2}\left(n\min
\biggl\{x^T Q x \ \bigg| \ x \in \mathbb{R}_+^{6!} ,
e^T x = 1\biggr\} - 9\right).
\end{align*}
\end{lemma}
{\bf Remark } For obvious reasons ($Q$ is a $720\times 720$--matrix)
we do not  include in this paper the matrix $Q$ in table form.  As we mentioned
above, $Q(\rho,\rho) = 9$ for every $\rho \in \Pi$, and therefore all the
diagonal entries of $Q$ are $9$.  It is not difficult to show that
$Q(\rho,\rho') \le 8$ if $\rho \neq \rho'$, so every non--diagonal entry of $Q$
is at most $8$.
The calculation of the entries of $Q$, using the
definition of $Q(\cdot,\cdot)$, and taking its symmetries into account
(see Section~\ref{subsec:symms} on this)
takes only a few seconds of computer time.
%{\tt http://www./matrix.txt}.\marginri{Exact web site?}
\end{section}

\begin{section}{Finding a lower bound for the optimization
problem}\label{quadraticprogramming}
Our aim in this section is to find a (reasonably good) lower bound for the
quadratic programming problem with $m=7$
given in Lemma~\ref{7ninequality}, in order to
obtain a (reasonably good) lower bound for $\Cr(\ksn)$.  The main result in
this section is the following.

\begin{proposition}\label{thebound}
Let $Q$ be the $6!\times 6!$ matrix of the quadratic form
$Q(\cdot,\cdot)$.  Then
\begin{equation*}
\min \biggl\{x^T Q x \ \bigg| \ x \in \mathbb{R}_+^{6!} , e^T x = 1\biggr\}
\ge 4.3593.
\end{equation*}
\end{proposition}
We devote this section to the proof of Proposition~\ref{thebound}.
It involves computer calculations; more details on this are given
in Section~\ref{subsec:compres}.

\begin{subsection}{The standard quadratic programming problem}
The problem we have formulated is known as {\em standard quadratic
optimization problem}.  The standard quadratic optimization problem (standard
QP) is to find the global minimizers of
a quadratic form over the standard simplex, {\em i.e.}\ we consider
the global optimization problem
\begin{equation}
\label{sqp}
\underline{p} :=\min_{x      \in \Delta} x^T   Qx
\end{equation}
where $Q$ is an arbitrary symmetric $d\times d$ matrix, $e$ is the all ones
vector, and $\Delta$ is the standard simplex
in  $\LR^d$,
$$\Delta = \{ x     \in \LR^d_+ : e^T  x      = 1\}.$$

We will now reformulate the standard QP problem as a convex optimization
problem in conic form.  First, we will review the relevant convex cones as well
as the duality theory of conic optimization.
We define the following convex cones.
\begin{itemize}
\item The $d\times d$ symmetric matrices: \\ $ \cS_d =\left\{X \in
\LR^d\times\LR^d, \; X = X^T\right\} $; \item The $d\times d$
symmetric positive semidefinite matrices: \\ $ \cS_d^+ =\left\{X
\in  \cS_d, \; y^TXy \ge 0 \;\; \forall y \in \LR^d\right\} $;
\item The $d\times d$ symmetric copositive matrices: \\ $\cC_d
=\left\{X \in  \cS_d,   \; y^TXy \ge 0 \;\; \forall y \in \LR^d,
\; y \ge 0 \right\} $; \item The $d\times d$ symmetric completely
positive matrices:\\ $\cC_d^* =\left\{ X = \sum_{i=1}^k y_iy_i^T,
\; y_i \in \LR^d, \; y_i \ge 0 \; (i=1,\ldots,k) \right\} $; \item
The $d\times d$ symmetric nonnegative matrices:\\ $\cN_d =\left\{X
\in  \cS_d, \;  X_{ij} \ge 0 \;\; (i,j = 1,\ldots,d)  \right\} $.
\end{itemize}
Recall that the completely positive cone is the dual of the
copositive cone, and that the nonnegative and semidefinite cones
are self-dual for the inner product $\langle X,Y\rangle:=\Tr(XY)$,
where `$\Tr$' denotes the trace operator.

For a given cone $\cK_d$ and its dual cone $\cK_d^*$ we define the
primal and dual pair of conic linear programs:
\begin{align}
\tag{P} p^* & := \inf_{X\in \cK_d} \left\{ \Tr(CX)  \mid  \Tr(A_iX) = b_i \;
(i=1,\ldots,M) \right\}\\
\tag{D} d^* & := \sup_{y \in \LR^m} \left\{ b^Ty  \mid  \sum_{i=1}^M
y_iA_i + S = C, \; S \in \cK_d^* \right\}.
\end{align}
If $\cK_d =  \cS_d^+$ we refer to semidefinite programming, if
$\cK_d = \cN_d$ to linear programming, and if $\cK_d = \cC_d$ to
copositive programming.

The well-known conic duality theorem (see e.g. Renegar~\cite{Ren})
gives the duality relations between $(P)$ and $(D)$.
\begin{theorem}[Conic duality theorem]
If there exists an interior feasible solution $X^0 \in
\mbox{int}(\cK_d)$ of $(P)$, and a feasible solution of $(D)$ then
$p^*=d^*$ and the supremum in $(D)$ is attained. Similarly, if
there exist feasible $y^0,S^0$ for $(D)$ where $S^0 \in
\mbox{int}(\cK_d^*)$, and a feasible solution of $(P)$, then $p^* =
d^*$ and the infimum in $(P)$ is attained.
\end{theorem}

Optimization over the cones $ \cS_d^+$ and $\cN_d$ can be done in
polynomial time (to compute an $\e$-optimal solution), but 
some NP-hard problems can be formulated as
copositive programs, see e.g. De Klerk and Pasechnik~\cite{DeKPas}.

\subsubsection{Convex reformulation of the standard QP}
We rewrite problem (\ref{sqp}) in the following way:
\[
\underline{p} :=\min_{x      \in \Delta} \Tr(Qxx^T).
\]
Now we define the cone of matrices
\[
\cK = \left\{ X \in \cS_d \;  : \; X = xx^T, \; x \ge 0\right\}.
\]
Note that the requirement $x \in \Delta$ corresponds to $X \in \cK$ with $\Tr\left(ee^TX\right) = 1$.

We arrive at the following reformulation of problem (\ref{sqp}):
\beq
\underline{p} = \min \left\{\Tr (QX) \; : \; \Tr \left(ee^TX\right) = 1, \; X \in \cK\right\}.
\label{sqp2}
\eeq
The last step is to replace the cone $\cK$ by its convex hull, which is simply the cone of completely positive
matrices, i.e.
\[
\mbox{conv}\left(\cK\right) = \cC_d^*
=\left\{ X = \sum_{i=1}^k y_iy_i^T, \; y_i \in \LR^n, \; y_i \ge 0 \;
(i=1,\ldots,k) \right\}.
\]
Replacing the feasible set by its convex hull does not change the optimal value
of problem (\ref{sqp2}),
since its objective function is linear. Thus we obtain the well-known convex
reformulation:
\beq
\label{sqp3}
\underline{p} = \min \left\{\Tr (QX) \mid
\Tr \left(ee^TX\right) = 1, \; X \in \cC_d^* \right\}.
\eeq
The dual problem takes the form:
\beq
\label{sqp4}
\underline{p} = \max \left\{ t \mid Q - tee^T \in \cC_d\right\},
\eeq
where $\cC_d$ is the cone of copositive matrices, as before.
Note that both problems have the same optimal value, in view of the conic duality theorem.
\end{subsection}

\begin{subsection}{Exploiting group symmetries}\label{subsec:symms}
We can reduce the number of variables in the optimization problems in
(\ref{sqp3},\ref{sqp4}) considerably
by exploiting the invariance properties of the quadratic function $x^TQx$.
This will also prove to be computationally necessary for the problems 
we intend to solve.

Consider the situation where the matrix $Q$ is invariant under the action of a
group $G$ of order $k=|G|$ of permutation
matrices $P\in G$, in the sense that
\[
Q = P^TQP \quad\mbox{for all}\quad P\in G.
\]
Then we have
\begin{eqnarray*}
\underline{p} &=& \min \left\{\Tr (QX)
\mid \Tr \left(ee^TX\right) = 1, \; X \in \cC_d^* \right\} \\
         &=& \min \left\{\Tr \left(P^TQPX\right) \mid
          \Tr \left(Pee^TPX\right) = 1, \;
          X \in \cC_d^* \right\} \mbox{ for any $P\in G$ } \\
         &=& \min \left\{\Tr \left(QP^TXP\right) \mid
          \Tr \left(ee^TP^TXP\right) = 1, \;
          X \in \cC_d^* \right\} \mbox{ for any $P \in G$ } \\
         &=& \min \left\{\Tr \left(Q\frac{1}{k}\left[\sum_{P\in G}
         P^TXP\right]\right) \mid
         \Tr \left(ee^T\left[\frac{1}{k}\sum_{P\in G}P^TXP\right]\right) = 1,
              \; X \in \cC_d^* \right\}. \\
\end{eqnarray*}
We can therefore restrict the optimization to the subset of the
feasible set obtained by replacing each feasible $X$ by the {\em
group average} $\frac{1}{k}\sum_{P\in G} P^TXP$, i.e.\ replacing
$X$ by its image under what is known in invariant theory as the
Reynolds operator. Note that if $X \in \cC_d^*$, then so is its
image under the group average.

In particular, we wish to compute a basis for the so-called {\em fixed point
subspace}:
\[
\cA := \left\{Y \in \cS_d \mid
Y = \frac{1}{k}\sum_{P\in G} P^TXP, \; X \in \cS_d\right\}.
\]
Note that $Q$ and $ee^T$ are elements of $\cA$ (set $X=Q$, 
respectively, $X=ee^T$).
Hence $Q-tee^T\in\cA$ for any $t$, and
\[
\underline{p} = \max \left\{ t \mid Q - tee^T \in \cC_d\right\}
              = \max \left\{  t \mid Q - tee^T \in \cC_d \cap \cA \right\}. 
\]
The right-hand side here is the dual of the primal 
problem when it was restricted to $\cA$ as above.

The next step is to compute a basis for the subspace $\cA$.
\end{subsection}

\begin{subsection}{Computing a basis for the fixed point subspace}
We assume for simplicity that $G$ acts transitively as a permutation group on
the standard basis vectors. (This holds in our setting.
A more general, and computationally less efficient, setting, can be found 
in Gatermann and Parrilo~\cite{GaPa}.)
The theory here is well-known, and goes back to Burnside, Schur and Wielandt.
See e.g.\ Cameron \cite{Cameron} for details.
Although we need a basis of $\cA$, the subspace of {\em symmetric} matrices
fixed by $G$, it is more natural to compute the basis $\mathcal{X}$
of the subspace $\cBB$ of {\em all}
fixed by $G$ matrices, and then pass on to $\cA$.

The dimension of $\cBB$ equals the number $r$ of orbits 
of $G$ on the Cartesian square of the standard basis.
The set of the latter orbits, also known as 2-orbits, 
naturally corresponds to certain set $\mathcal{X}$ of 
$d\times d$ zero-one matrices. Namely, for each $X\in\mathcal{X}$ one has
$X_{ij}=1$ if and only if  $X_{P(i),P(j)}=1$ for all $P\in G$
and all $1\leq i\leq j\leq |\Pi|$.
As $G$ is transitive on the standard basis vectors,
the identity matrix $I$ belongs to $\mathcal{X}$.
We also have $\sum_{X\in\mathcal{X}} X=ee^T$.

As $\mathcal{X}$  is closed under the matrix transposition,  
i.e. $X^T\in\mathcal{X}$ for any $X\in\mathcal{X}$,
$${\mathcal X}_{\cA}=\{A_1,\dots,A_M\}=
\{X\mid X=X^T\in{\mathcal X}\}\cup 
\{X+X^T\mid X\in{\mathcal X},\ X\ne X^T\}$$ 
is a basis of $\cA$.  Each $A\in\mathcal{X}_{\cA}$ 
is a symmetric zero-one matrix and 
$\sum_{A\in{\mathcal X}_{\cA}}A = ee^T$.
Moreover
\[
\left\{Y \in \cS_d \mid  Y = \sum_{i=1}^M y_iA_i\right\} = \cA \equiv
\left\{Y \in \cS_d \mid  Y = \frac{1}{k} \sum_{P\in G} P^TXP, \;
X \in \cS_d\right\}.
\]
Since $Q \in \cA$, we will write $Q = \sum_{i=1}^M b_iA_i$.

It is worth mentioning that algebraically the vector space $\cBB$ behaves very
nicely:
it is closed under multiplication. In other words, $\cBB$ is a matrix
algebra of dimension $r$, also known as the {\em centralizer ring}
of the permutation group $G$.

We proceed to describe $G$ and $\cBB$ in our case.
For us $G$ is isomorphic to the direct product
$Sym(m)\times Sym(2)$, of symmetric
groups $Sym(m)$ and $Sym(2)$,
where $Sym(m)$ acts (as a permutation group) by conjugation on the $d=(m-1)!$
elements of %the conjugacy class 
$\Pi$, %of the $m$-cycle $(0,1,\dots,m-1)$, 
and
$Sym(2)$ acts (as a permutation group) on $\Pi$ by switching
$\pi\in\Pi$ with $\pi^{-1}\in\Pi$.

Computing $\mathcal{X}$ is an elementary combinatorial procedure,
that can be found in one form  or another in many computer algebra
systems, so one does not have to program this again.
First, the permutations that generate $Sym(m)\times Sym(2)$ in its action
on $\Pi$ are computed. The action of $Sym(2)$ is already known, and
is described by the permutation $g_0$, say.
In its usual action on $m$ symbols, $Sym(m)$ is generated by
$h_1=(0,1,\dots,m-1)$ and $h_2=(0,1)$. 
These $h_i$ (for $i=1,2$) act on $\Pi$ by 
mapping each $\pi\in\Pi$ to
$h_i\pi h_i^{-1}$. Denote by $g_i$ (for $i=1,2$) the
permutations of $\Pi$ that realize these actions.

Next,  one computes the orbits of the permutation group
$Sym(m)\times Sym(2)=\langle g_0, g_1, g_2\rangle$ on the Cartesian
square $\Pi\times\Pi$ of $\Pi$, by ``spinning''
$(\pi_i,\pi_j)\in\Pi\times\Pi$: begin with
$S_{ij}=\{(\pi_i,\pi_j)\}$ and apply the generators
$g_i$, $0\leq i\leq 2$, in a loop until $S_{ij}$ stops growing.
Then one sets $\Pi:=\Pi-S_{ij}$ and repeats until $\Pi$ is exhausted.

\medskip

When $m=7$ one has $r=78$ and $M=56$. Note that here the algebra $\cBB$ is not
commutative.

When $m=5$ one has $r=M=6$, and $\cBB$ commutative.
\end{subsection}

\begin{subsection}{Reformulation of the optimization problem}
We can now reformulate the dual problem by using the basis of $\cA$ to obtain:
\[
\underline{p} = 
\max \left\{  t \mid Q - tee^T \in \cC_d \cap \cA \right\} 
              = \max \left\{  t \mid \sum_{i=1}^M (b_i-t)A_i \in
          \cC_d \right\}.
\]
We will now proceed to derive a lower bound on $\underline{p}$ by solving the
dual problem approximately.
\end{subsection}

\newcommand{\MM}{S}
\newcommand{\rr}{\ell}

\begin{subsection}{Approximations of the copositive cone}
The problem of determining whether a matrix is not copositive is
NP-complete, as shown by Murty and Kabadi~\cite{Murty87}. 
We therefore wish to replace the copositive cone
$\cC_d$ by a conic subset, in such a way that the resulting optimization
problem becomes tractable.
We can represent the copositivity requirement for a $(d\times d)$
symmetric matrix $\MM$ as
\beq \label{copositivity} P(x) :=(x \circ
x)^T \MM(x \circ x) = \sum_{i,j=1}^d \MM_{ij}x_i^2x_j^2 \ge 0, \;
\forall x \in \LR^d,
\eeq
where `$\circ$' indicates the
componentwise (Hadamard) product.
We therefore wish to know whether the polynomial $P(x)$ is
nonnegative for all $x \in \LR^d$.
Although one apparently
cannot answer this question in polynomial time in general,
as it is an NP-hard problem,
one can decide using semidefinite programming whether $P(x)$ can be written
as a sum of squares.

Parrilo \cite{Parrilo2000} showed that $P(x)$ in \reff{copositivity} allows a
sum of squares decomposition if and only if $\MM \in  \cS_d^+ + \cN_d$,
which is a well-known sufficient condition for copositivity.
Set $\cK^0_d$ to be the convex cone $\cK^0_d  = \cS_d^+ + \cN_d$.

Higher order sufficient conditions can be derived by considering
the polynomial:
\beq \label{higher} P^{(\rr)}(x) =
P(x)\left(\sum_{i=1}^d x_i^2 \right)^{\rr} = \left(\sum_{i,j=1}^d
\MM_{ij}x_i^2x_j^2\right)\left(\sum_{i=1}^d x_i^2 \right)^{\rr}, \eeq
and asking
whether $P^{(\rr)}(x)$  -- which is a homogeneous polynomial of
degree $2(\rr+2)$ -- has a sum of squares decomposition, or whether
it only has nonnegative coefficients.

For $\rr=1$, Parrilo \cite{Parrilo2000} showed that a sum of squares
decomposition exists if and only
if\footnote{In fact, Parrilo \cite{Parrilo2000} only proved the
'if'-part; the converse is proven in
Bomze and De Klerk \cite{BomDeK}.} the following system of linear
matrix inequalities has a solution:
\begin{eqnarray}
\MM - \MM^{(i)} & \in &  \cS_d^+, \;\;\; i=1,\ldots,d \label{def:k1a1}\\
\MM^{(i)}_{ii} &=& 0, \;\;\; i=1,\ldots,d \label{def:k1a2}\\
\MM^{(i)}_{jj}+2\MM^{(j)}_{ij} &=& 0, \;\;\; i\neq j \label{def:k1a3}\\
\MM^{(i)}_{jk}+\MM^{(j)}_{ik}+\MM^{(k)}_{ij} &\ge& 0, \;\;\;
i< j < k, \label{def:k1a4}
\end{eqnarray}
where $\MM^{(i)}$ $(i=1,\ldots,d)$ are symmetric matrices.  Similar to the
$\rr=0$ case, we define $\cK^1_d$ as the (convex) cone of matrices $\MM$
for which the above system has a solution.

We will consider the lower bounds we
get by replacing the copositive cone by either $\cK_d^0$ or $\cK_d^1$:
\beq
\label{relax}
\underline{p} \ge p_{\rr} :=
\max \left\{ t \mid Q - tee^T \in \cK^{\rr}_d\right\},\qquad \rr \in \{0,1\}.
\eeq
\end{subsection}

\begin{subsection}{Approximations (relaxations) of the copositive cone}
We will now study the relaxation obtained by replacing the copositive cone by
its proper subset $\cK_d^0$.
In other words, we study the relaxation:
\begin{eqnarray*}
\underline p &=& \max \left\{  t \mid \sum_{i=1}^M (b_i-t)A_i \in \cC_d \right\} \\
             &\ge& p_0 :=
         \max \left\{  t \mid \sum_{i=1}^M (b_i-t)A_i \in \cK_d^0 = \cS_d^+ + \cN_d \right\}
\end{eqnarray*}
We rewrite $\sum_{i=1}^M (b_i-t)A_i \in \cK_d^0$ as
\[
\sum_{i=1}^M (b_i-t)A_i = \sum_{i=1}^M y_iA_i + \sum_{i=1}^M z_iA_i,
\quad\text{ where} \sum_{i=1}^M y_iA_i \in \cS_d^+
\text{ and } \sum_{i=1}^M z_iA_i \in \cN_d.
\]
Note that, since the $A_i$'s are 0-1 matrices that sum to $ee^T$, it follows that $z_i \ge 0$.
Moreover,
\begin{equation*}
b_i - t = y_i + z_i \qquad\text{implies}\qquad b_i - t - y_i \ge 0.
\end{equation*}
We obtain the relaxation:
\beq\label{K0_rel}
p_0 = \max \left\{  t \mid b_i - t - y_i \ge 0 \; (i=1,\ldots,M), \;
\sum_{i=1}^M y_iA_i \in \cS_d^+ \right\}.
\eeq
\end{subsection}

\begin{subsection}{Block factorization}
The next step in reducing the problem size is to perform a similarity
transformation that simultaneously block-diagonalizes the matrices
$A_1,\ldots,A_M$.
In particular, we want to find an orthogonal matrix $V$ 
%(a matrix is orthogonal if
%$V^{-1}=V^T$) 
such that the matrices
\[
\tilde A_i := VA_i V^{-1} \;\;\; i= 1,\ldots,M,
\]
all have the same block diagonal structure, 
and the maximum block size is as small as possible. 
Note that the conjugation preserves spectra, and
orthogonality of $V$ preserves symmetry.

This will further reduce the size of the relaxation (\ref{K0_rel}) via
\begin{eqnarray*}
p_0 &=& \max \left\{  t \mid b_i - t - y_i \ge 0 \;
(i=1,\ldots,M), \; \sum_{i=1}^M y_iA_i \in \cS_d^+ \right\} \\
  &=& \max \left\{  t \mid b_i - t - y_i \ge 0 \; (i=1,\ldots,M), \;
  \sum_{i=1}^M y_iVA_i V^{-1}  \in \cS_d^+ \right\} \\
 &=& \max \left\{  t \mid b_i - t - y_i \ge 0 \; (i=1,\ldots,M), \;
 \sum_{i=1}^M y_i\tilde A_i \in \cS_d^+ \right\}. \\
\end{eqnarray*}

The necessity to restrict to orthogonal $V$'s lies in the fact that there is
currently no software (or algorithms) available that would be able to deal with
non-symmetric $\tilde A_i$'s.

Computing the finest possible block decomposition (this would mean finding
explicitly the orthogonal bases for the irreducible submodules of the natural
module of $G$ in its action by the matrices $P$) is computationally not easy,
especially due to the orthogonality requirement on $V$. We restricted ourselves
to decomposing into two blocks of equal size $\frac{d}{2}\times\frac{d}{2}$.
Namely, each row corresponds to a cyclic permutation $g\in\Pi$, and the natural
pairing $(g,g^{-1})$ can be used to construct $V=\frac{\sqrt{2}}{2}V'$, as
follows:
\begin{itemize}
\item the first half of the rows of $V'$ are
characteristic vectors of the 2-subsets $\{g,g^{-1}\}$, $g\in\Pi$;
\item the second half of the rows of $V'$ consists of ``twisted''
rows from the first
half: namely one of the two 1s is replaced by -1.
\end{itemize}
It is obvious that $V'V'^T=2I$ and thus $V$ is orthogonal.

\paragraph{Remark}
It is worth mentioning that in \cite{Schr79} Schrijver essentially dealt,
in a different context, with a
similar setup, except that in his case the elements of the basis
$\mathcal{X}$ of $\cBB$ were symmetric
and (hence) the algebra $\cBB$ commutative.
In such a situation the elements of $\mathcal{X}$ 
can be simultaneously diagonalized, and the
corresponding optimization problem becomes a linear programming problem.
\end{subsection}

\begin{subsection}{Computational results: proof of Theorem~\ref{thebound}}
\label{subsec:compres}
%{\tt this section is just to relate the results we have obtained
%to the relevant relaxations.}
The combinatorial/group theoretic part of the computations,
namely of the $A_i$'s, $V$, and $Q=\sum_i b_i A_i$ was
performed using a computer algebra system {\sf GAP} \cite{GAP}
version 4.3, and its shared package {\sf GRAPE} by Soicher~\cite{GRAPE}.
Semidefinite programs were solved using {\sf SeDuMi} by 
Sturm~\cite{sedumi}, version 1.05 under {\sf Matlab}~6.5.
The biggest SDP took about 10 minutes of CPU time of a Pentium~4
with 1~GB of RAM.

In addition, the results were verified using Maple.
Namely, for $t=p_0$ and $y$, the variables computed upon
solving (\ref{K0_rel}), we checked that the corresponding
(matrix and scalar) inequalities in (\ref{K0_rel}) hold.
As $p_0$ is a lower bound on $\underline{p}$, we thus validated
the computed value of $p_0$ independently of
the SDP solver used.

For the test case of $K_{5,n}$ we solved the relaxed problem
(\ref{relax}) with $\rr=1$ to obtain% the optimal value
\[
p_1 \approx  1.9544,\quad\text{that is}\quad
\mbox{cr}(K_{5,n}) \ge \hlf(1.9544)n^2 = 0.9772n^2,
\]
asymptotically.
The correct asymptotic value is known to be $\mbox{cr}(K_{5,n}) = n^2$, which
shows the quality of the bound.
The weaker bound for $\rr=0$ in (\ref{relax}) yields,
still quite tight
\[
p_1 \approx  1.94721,\quad\text{that is}\quad
\mbox{cr}(K_{5,n}) \ge \hlf(1.94721)n^2 = 0.973605n^2.
\]
For the case $K_{7,n}$ we solved the relaxed problem (\ref{relax}) with
$\rr=0$ to obtain% the optimal value
\[
p_0 \approx 4.3593,\quad\text{that is}\quad
\mbox{cr}(K_{7,n}) \ge \hlf(4.3593)n^2 = 2.1796n^2,
\]
asymptotically.

%{\tt The relaxation with $\rr=1$ is too
%large too solve in this case, but we can solve a
%relaxation of this relaxation. But this will have to wait ... }

\paragraph{Proof of Theorem~\ref{maintheorem}.}
For the sake of completeness, we close this section with the observation that
Theorem~\ref{maintheorem} has been proved.
It follows from
Lemmas~\ref{7ninequality} and~\ref{thebound}. \qed
\end{subsection}
\end{section}

\begin{section}{Improved bounds for the crossing numbers of $K_{m,n}$ and
$K_n$}\label{consequences}
Perhaps the most appealing consequence of our improved bound for $\Cr(\ksn)$ is
that it also allows us to give improved lower bounds for the crossing numbers
of $K_{m,n}$ and $K_n$. The quality of the new bounds is perhaps
best appreciated in terms of the following asymptotic parameters:
\begin{equation*}
A(m):= \lim_{n\to\infty}\frac{\Cr(K_{m,n})} {Z(m,n)}, \qquad\qquad
B:= \lim_{n\to\infty}\frac{\Cr(K_{n,n})}{Z(n,n)},
\end{equation*}
(see Richter and Thomassen~\cite{richterthomassen}).
These natural parameters give us a good idea of our current standing with
respect to Zarankiewicz's Conjecture.  It is not difficult to show that $A(m)$
(for every integer $m \ge 3$) and $B$ both exist \cite{richterthomassen}.

Previous to the new bound we report in Theorem~\ref{maintheorem}, the best
known lower bounds for $A(m)$ and $B$ were $A(m) \ge  0.8\frac{m}{m-1}$ and
(consequently) $B \ge 0.8$.  Both bounds were obtained by using the known value
of $\Cr(K_{5,n})$, and applying a standard counting argument.

By applying the same counting argument, but instead using the bound given by
Theorem~\ref{maintheorem}, we improve these asymptotic quotients to
$A(m) > 0.83\frac{m}{m-1}$ and $B > 0.83.$

The improved lower bound for $B$ has an additional, important application.  It
has been long conjectured that $\Cr(K_n) =Z(n)$, where
\begin{equation*}
Z(n)=\frac{1}{4}\bigfloor{\frac{n}{2}}\bigfloor{\frac{n-1}{2}}
\bigfloor{\frac{n-2}{2}}\bigfloor{\frac{n-3}{2}},
\end{equation*}
but this has been verified only for $n \le 10$  (see for
instance~\cite{erdosguy}).   As we did with $K_{m,n}$, it is natural to inquire
about the asymptotic parameter
\begin{equation*}
C:= \lim_{n\to\infty} \frac{\Cr(K_n)}{Z(n)}.
\end{equation*}

In~\cite{richterthomassen} it is proved that $C$ exists, and, moreover, $C \ge
B$.  In view of this, our improved lower bound for $B$ yields
$C > 0.83.$

We summarize these results in the following statement.
\begin{theorem}
With $Z(m,n)$ and $Z(n)$ as above,
%\begin{align*}
$$
\lim_{n\to\infty} \frac{\Cr(K_{m,n})}{Z(m,n)}  \ge  0.83\frac{m}{m-1},\quad
   \lim_{n\to\infty} \frac{\Cr(K_{n,n})}{Z(n,n)} \ge  0.83,
   \quad \text{\rm and }
 \lim_{n\to\infty}\frac{\Cr(K_n)}{Z(n)} \ge     0.83  \eqno\square
 $$
% \end{align*}
 \end{theorem}

We close this section with a few words on some important recent developments involving the {\it rectilinear} crossing number of $K_n$.

The {\it rectilinear crossing number} $\rcr(G)$ of a graph $G$ is the minimum number of
pairwise intersections of edges in a drawing of $G$ in the plane, with the additional restriction that all edges of $G$ must be drawn as straight segments.

It is known that $\rcr(K_n)$ and $\Cr(K_n)$ may be different (for instance, $\rcr(K_8) = 19$, whereas $\Cr(K_8) = 18$).
  While we have a (non--rectilinear) way of drawing $K_n$ that shows $\Cr(K_n) \le Z(n)$ (equality is conjectured to hold, as we observed above), good upper bounds for $\rcr(K_n)$ are notoriously difficult to obtain.  Currently, the best upper bound known is $\rcr(K_n) \le 0.3807{n\choose 4}$ 
  (see Aichholzer {\em et al.}~\cite{aurenhammeretal}).

  For many years the best lower bounds known for $\rcr(K_n)$ were considerably smaller (around $0.32{n\choose 4}$) than the best upper bounds available 
(currently around $0.380{n\choose 4}$). 
However, remarkably better lower bounds have been recently proved  
independently by \'Abrego and Fern\'andez--Merchant~\cite{abrego} and by 
Lov\'asz {\em et al.}~\cite{lovasz}, and
refined by Balogh and Salazar~\cite{bs}.  In~\cite{abrego}, the   technique
of allowable sequences is used to show that $\rcr(K_n) \ge 0.375{n\choose 4}$.
Lov\'asz {\em et al.} use similar methods to prove 
$\rcr(K_n) > 0.37501{n\choose 4} + O(n^3)$.  
Recently, Balogh and Salazar  
improved this to $\rcr(K_n) > 0.37553{n\choose 4} + O(n^3)$~\cite{bs}.  
The importance of establishing that
$\rcr(K_n)$ is strictly greater than $0.375{n\choose 4} + O(n^3)$ 
is that it effectively shows that the ordinary and the 
rectilinear crossing numbers of $K_n$ are different in the 
asymptotically relevant term, namely $n^4$.

\end{section}

\paragraph{Acknowledgements.}
Etienne de Klerk would like to thank Pablo Parrilo for his valuable comments.

 \end{document}